\documentclass{amsart}
\usepackage{amssymb}
\usepackage{amsfonts}
\usepackage{latexsym}

\newcommand{\C}{{\mathbb C}}

\newcommand{\Gr}{\operatorname{Gr}}

\newtheorem{thm}{Theorem}
\newtheorem{cor}{Corollary}
\newtheorem{lemma}{Lemma}

\theoremstyle{definition}

\newcommand{\refthm}[1]{Theorem~\ref{#1}}

\begin{document}

\title[Eigenvalues of matrices with sum of bounded rank]
{Eigenvalues of Hermitian matrices with \\ positive sum of bounded rank}

\author{Anders Skovsted Buch}
\address{Matematisk Institut\\
         Aarhus Universitet\\
         Ny Munkegade, 8000 {\AA}rhus C, Denmark}
\email{abuch@imf.au.dk}

\begin{abstract}
We give a minimal list of inequalities characterizing the possible
eigenvalues of a set of Hermitian matrices with positive semidefinite
sum of bounded rank.  This answers a question of A.~Barvinok.
\end{abstract}

\date{September 1, 2004.}
\maketitle

\newcommand{\decr}{$(\dagger)$}
\newcommand{\major}[1]{$(\triangleright_{#1})$}
\newcommand{\rkcond}[2]{$(\triangleleft_{#1,#2})$}

\section{Introduction}

The combined work of A.~Klyachko \cite{klyachko:stable}, A.~Knutson,
T.~Tao \cite{knutson.tao:honeycomb} and C.~Woodward
\cite{knutson.tao.ea:honeycomb}, and P.~Belkale \cite{belkale:local}
produced a minimal list of inequalities determining when three
(weakly) decreasing $n$-tuples of real numbers can be the eigenvalues
of Hermitian $n\times n$ matrices which add up to zero.  We refer to
\cite{fulton:eigenvalues*1} for a description of this work, as well as
references to earlier work and applications to a surprising number of
other mathematical disciplines.

S.~Friedland applied these results to determine when three decreasing
$n$-tuples of real numbers can be the eigenvalues of three Hermitian
matrices with positive semidefinite sum, that is, the sum should have
non-negative eigenvalues \cite{friedland:finite}.  Friedland's answer
included the inequalities of the above named authors, except that a
trace equality was changed to an inequality.  Friedland's result also
needed some extra inequalities.  W.~Fulton has proved
\cite{fulton:eigenvalues*2} that the extra inequalities are
superfluous, and that the remaining ones form a minimal list, i.e.\ 
they correspond to the facets of the cone of permissible eigenvalues.
All of these results have natural generalizations that work for any
number of matrices \cite{fulton:eigenvalues*1,
  knutson.tao.ea:honeycomb}.

In this paper we address the following more general question, which
was formulated by A.~Barvinok and passed along to us by Fulton.  Given
weakly decreasing $n$-tuples of real numbers $\alpha(1), \dots,
\alpha(m)$ and an integer $r \leq n$, when can one find Hermitian
$n\times n$ matrices $A(1), \dots, A(m)$ such that $\alpha(s)$ is the
eigenvalues of $A(s)$ for each $s$ and the sum $A(1)+\dots+A(m)$ is
positive semidefinite of rank at most $r$?  The above described
problems correspond to the extreme cases $r=0$ and $r=n$.

Let $\alpha(1), \alpha(2), \dots, \alpha(m)$ be $n$-tuples of reals,
with $\alpha(s) = (\alpha_1(s), \dots, \alpha_n(s))$.  The requirement
that these $n$-tuples should be decreasing is equivalent to the
inequalities
\begin{equation} \tag{\text{$\dagger$}}
  \alpha_1(s) \geq \alpha_2(s) \geq \dots \geq \alpha_n(s)
\end{equation}
for all $1 \leq s \leq m$.

Given a set $I = \{a_1 < a_2 < \dots < a_t\}$ of positive integers, we
let $s_I = \det(h_{a_i-j})_{t\times t}$ be the Schur function for the
partition $\lambda(I) = (a_t-t, \dots, a_2-2, a_1-1)$.  Here $h_i$
denotes the complete symmetric function of degree $i$.  Fulton's
result \cite{fulton:eigenvalues*2} states that the $n$-tuples
$\alpha(1),\dots,\alpha(m)$ can be the eigenvalues of Hermitian
matrices with positive semidefinite sum if and only if
\begin{equation} \tag{\text{$\triangleright_n$}}
  \sum_{s=1}^m \sum_{i \in I(s)} \alpha_{i}(s) \geq 0
\end{equation}
for all sequences $(I(1),\dots,I(m))$ of subsets of $[n] =
\{1,2,\dots,n\}$ of the same cardinality $t$ ($1 \leq t \leq n$), such
that the coefficient of $s_{\{n-t+1,n-t+2,\dots,n\}}$ in the Schur
expansion of the product $s_{I(1)} s_{I(2)} \cdots s_{I(m)}$ is equal
to one.  Notice that this coefficient is one if and only if the
corresponding product of Schubert classes on the Grassmannian
$\Gr(t,\C^n)$ equals a point class.

The added condition that the rank of the sum of matrices is at most
$r$ results in the additional inequalities
\begin{equation} \tag{\text{$\triangleleft_{n,r}$}}
  \sum_{s=1}^m \sum_{p \in P(s)} \alpha_{n+1-p}(s) \leq 0
\end{equation}
for all sequences $(P(1),\dots,P(m))$ of subsets of $[n-r]$ of the
same cardinality $t$ ($1 \leq t \leq n-r$), such that
$s_{\{n-r-t+1,\dots,n-r\}}$ has coefficient one in the product
$s_{P(1)} s_{P(2)} \cdots s_{P(m)}$.  Equivalently, a product of
Schubert classes on $\Gr(t,\C^{n-r})$ should be a point class.  We
remark that without the requirement that a Hermitian matrix is
positive semidefinite, rank conditions on the matrix do not correspond
to linear inequalities in the eigenvalues.  The following theorem is
our main result.

\begin{thm} \label{T:majorrank}
  Let $\alpha(1),\dots,\alpha(m)$ be $n$-tuples of real numbers
  satisfying \decr{}, and let $r \leq n$ be an integer.  There exist
  Hermitian $n \times n$ matrices $A(1), \dots, A(m)$ with eigenvalues
  $\alpha(1), \dots, \alpha(m)$ such that the sum $A(1)+\cdots+A(m)$
  is positive semidefinite of rank at most $r$, if and only if the
  inequalities \major{n} and \rkcond{n}{r} are satisfied.
  Furthermore, for $r \geq 1$ and $m \geq 3$ the inequalities \decr{},
  \major{n}, and \rkcond{n}{r} are independent in the sense that they
  correspond to facets of the cone of admissible eigenvalues.
\end{thm}

As proved in \cite{knutson.tao.ea:honeycomb}, the minimal set of
inequalities in the case $r=0, m\geq 3$ consists of the inequalities
\major{n} for $t<n$, along with the trace equality $\sum_{s=1}^m
\sum_{i=1}^n \alpha_i(s) = 0$ and, for $n > 2$, also the inequalities
\decr{}.  The cases $r=0, m\leq 2$, or $m=1$ are not interesting.  The
situation for $m=2$ and $r>0$ is described by the following special
cases of Weyl's inequalities \cite{weyl:das} (see also
\cite[p.~3]{fulton:eigenvalues*1}).

\begin{cor} \label{C:m2}
  Let $\alpha(1), \alpha(2)$ be $n$-tuples satisfying \decr{}, and let
  $r \leq n$ be an integer.  There exist Hermitian $n \times n$
  matrices $A(1), A(2)$ with eigenvalues $\alpha(1), \alpha(2)$ such
  that the sum $A(1)+A(2)$ is positive semidefinite of rank at most
  $r$, if and only if $\alpha_i(1) + \alpha_{j}(2) \geq 0$ for
  $i+j=n+1$ and $\alpha_i(1) + \alpha_{j}(2) \leq 0$ for $i+j=n+r+1$.
  These inequalities are independent when $r \geq 1$; they imply
  \decr{} for $r=1$, and are independent of \decr{} for $r \geq 2$.
\end{cor}
\begin{proof}
  Given subsets $I, J \subset [n]$ of cardinality $t$, the coefficient
  of $s_{\{n-t+1,\dots,n\}}$ in $s_I \cdot s_J$ is equal to one if and
  only if $J = \{ n+1-i \mid i \in I \}$.  This implies that the
  inequalities \major{n} and \rkcond{n}{r} are consequences of the
  inequalities of the corollary.  The claims about independence of
  inequalities are left as an easy exercise.
\end{proof}

Theorem~\ref{T:majorrank} also has the following consequence.
Although the statement does not use any inequalities, it appears to be
non-trivial to prove without the use of inequalities.

\begin{cor}
  Let $\alpha(1),\dots,\alpha(m)$ be $n$-tuples of real numbers and
  let $r \leq n$.  There exist Hermitian $n \times n$ matrices
  $A(1),\dots,A(m)$ with these eigenvalues such that $A(1)+\dots+A(m)$
  is positive semidefinite of rank at most $r$, if and only if there
  are Hermitian $n\times n$ matrices with the same eigenvalues and
  positive semidefinite sum, as well as Hermitian $(n-r) \times (n-r)$
  matrices $C(1),\dots,C(m)$ with negative semidefinite sum, such
  that the eigenvalues of $C(s)$ are the $n-r$ smallest numbers from
  $\alpha(s)$.
\end{cor}
\begin{proof}
  The inequalities \rkcond{n}{r} for $n$-tuples
  $\alpha(1),\dots,\alpha(m)$ are identical to the inequalities
  \major{n-r} for $\tilde \alpha(1),\dots,\tilde \alpha(m)$, where
  $\tilde \alpha(s) = (-\alpha_n(s) \geq \dots \geq
  -\alpha_{r+1}(s))$.
\end{proof}

Our proof of \refthm{T:majorrank} is by induction on $r$, where we
rely on the above mentioned results of Klyachko, Knutson, Tao,
Woodward, and Belkale to cover the base case $r=0$.  To carry out the
induction, we use an enhancement of Fulton's methods from
\cite{fulton:eigenvalues*2}.  We remark that \refthm{T:majorrank}
remains true if the Hermitian matrices are replaced with real
symmetric matrices or even quaternionic Hermitian matrices.  This
follows because the results for zero-sum matrices hold in this
generality \cite[Thm.~20]{fulton:eigenvalues*1}.

We thank Barvinok and Fulton for the communication of Barvinok's
question, and Fulton for many helpful comments to our paper.

\section{The inequalities are necessary and sufficient}

In this section we prove that the inequalities of
Theorem~\ref{T:majorrank} are necessary and sufficient.
For a subset $I = \{ a_1 < a_2 < \dots < a_t \}$ of $[n]$ of
cardinality $t$, we let $\sigma_I \in H^* \Gr(t,\C^n)$ denote the
Schubert class for the partition $\lambda(I) = (a_t-t,\dots,a_1-1)$.
The corresponding Schubert variety is the closure of the subset of
points $V \in \Gr(t,\C^n)$ for which $V \cap \C^{n-a_i} \subsetneq V
\cap \C^{n-a_i+1}$ for all $1 \leq i \leq t$.  Let $S^n_t(m)$ denote
the set of sequences $(I(1),\dots,I(m))$ of subsets of $[n]$ of
cardinality $t$, such that the product $\prod_{s=1}^n \sigma_{I(s)}$
is non-zero in $H^* \Gr(t,\C^n)$, and we let $R^n_t(m) \subset
S^n_t(m)$ be the subset of sequences such that $\prod_{s=1}^n
\sigma_{I(s)}$ equals the point class
$\sigma_{\{n-t+1,\dots,n-1,n\}}$.

The inequalities \major{n} are indexed by all sequences
$(I(1),\dots,I(m))$ which belong to the set $R^n(m) = \bigcup_{1\leq
  t\leq n} R^n_t(m)$.  Furthermore, it is known \cite{belkale:local,
  knutson.tao.ea:honeycomb} that if $\alpha(1),\dots,\alpha(m)$ are
decreasing $n$-tuples of reals satisfying \major{n}, then they also
satisfy the larger set of inequalities indexed by sequences from
$S^n(m) = \bigcup_{1\leq t\leq n} S^n_t(m)$, that is $\sum_{s=1}^m
\sum_{i\in I(s)} \alpha_i(s) \geq 0$ for all $(I(1),\dots,I(m)) \in
S^n(m)$.  Similarly, the inequalities of \rkcond{n}{r} are indexed by
$R^{n-r}(m)$, and if $\alpha(1),\dots,\alpha(m)$ satisfy these
inequalities, then we also have $\sum_{s=1}^m \sum_{p\in P(s)}
\alpha_{n+1-p}(s) \leq 0$ for all sequences $(P(1),\dots,P(m)) \in
S^{n-r}(m)$.

We first show that the inequalities \major{n} and \rkcond{n}{r} are
necessary.  Suppose $A(1),\dots,A(m)$ are Hermitian $n \times n$
matrices with eigenvalues $\alpha(1),\dots,\alpha(m)$, such that the
sum $B = A(1)+\dots+A(m)$ is positive semidefinite with rank at most
$r$.  Let $\beta = (\beta_1 \geq \dots \geq \beta_r,0,\dots,0)$ be the
eigenvalues of $B$.  For any sequence $(I(1),\dots,I(m)) \in R^n_t(m)$
we have that $(J,I(1),\dots,I(m))$ is in $R^n_t(m+1)$ where $J =
\{1,2,\dots,t\}$.  This is true because $\sigma_J \in H^* \Gr(t,\C^n)$
is the unit.  Since $-B+A(1)+\dots+A(m) = 0$, it follows from
\cite[Thm.~11]{fulton:eigenvalues*1} that
\[ - \sum_{j \in J} \beta_{n+1-j} + \sum_{s=1}^m \sum_{i \in I(s)}
   \alpha_i(s) \geq 0 \,,
\]
which implies \major{n} because each $\beta_j$ is non-negative.

On the other hand, if $(P(1),\dots,P(m)) \in R^{n-r}_t(m)$, then
$(Q,P(1),\dots,P(m)) \in R^n_t(m)$ where $Q = \{ r+1,r+2,\dots,r+t
\}$.  This follows from the Littlewood-Richardson rule, since
$\lambda(Q) = (r)^t$ is a rectangular partition with $t$ rows and $r$
columns.  Since $B - A(1) - \dots - A(m) = 0$,
\cite[Thm.~11]{fulton:eigenvalues*1} implies that
\[ \sum_{q \in Q} \beta_q - \sum_{s=1}^m \sum_{p \in P(s)}
   \alpha_{n+1-p}(s) \geq 0 \,.
\]
Since $\beta_q = 0$ for every $q \in Q$, this shows that \rkcond{n}{r}
is true.

If $I = \{ i_1 < i_2 < \dots < i_t \}$ is a subset of $[n]$ and $P$ is
a subset of $[t]$, we set $I_P = \{ i_p \mid p \in P \}$.  To prove
that the inequalities are sufficient, we need the following
generalization of \cite[Prop.~1 (i)]{fulton:eigenvalues*2}.

\begin{lemma} \label{L:index}
  Let $(I(1),\dots,I(m)) \in S^n_t(m)$ and let $(P(1),\dots,P(m)) \in
  S^{t-r}_x(m)$, where $0 \leq r \leq t$.  Then $(I(1)_{P(1)}, \dots,
  I(m)_{P(m)})$ belongs to $S^{n-r}_x(m)$.
\end{lemma}
\begin{proof}
  The case $r=0$ of this Lemma is equivalent to part (i) of
  \cite[Prop.~1]{fulton:eigenvalues*2}.  We deduce the lemma from this
  case using straightforward consequences of the Littlewood-Richardson
  rule.
  
  Set $Q = \{p+r \mid p \in P(1)\}$.  Since $\lambda(Q) = (r)^x +
  \lambda(P(1))$, it follows that $\sigma_Q \cdot \prod_{s=2}^m
  \sigma_{P(s)} \neq 0$ on $\Gr(x,t)$.  By the $r=0$ case, this
  implies that $\sigma_{I(1)_Q} \cdot \prod_{s=2}^m
  \sigma_{I(s)_{P(s)}} \neq 0$ on $\Gr(x,n)$.  Now notice that if
  $P(1) = \{p_1<\dots<p_x\}$ and $I(1) = \{i_1 < \dots < i_t\}$ then
  the $j$th element of $I(1)_{Q}$ is $i_{p_j+r} \geq i_{p_j}+r$, i.e.\ 
  $\lambda(I(1)_Q) \supset (r)^x + \lambda(I(1)_{P(1)})$.  This means
  that $\sigma_{(r)^x+\lambda(I(1)_{P(1)})} \cdot \prod_{s=2}^m
  \sigma_{I(s)_{P(s)}}$ is also non-zero on $\Gr(x,n)$, which implies
  that $\prod_{s=1}^m \sigma_{I(s)_{P(s)}} \neq 0$ on $\Gr(x,n-r)$.
\end{proof}

We also need the following special case of Corollary~\ref{C:m2},
which comes from reformulating the Pieri rule in terms of eigenvalues.

\begin{lemma} \label{L:rankone}
Let $\alpha = (\alpha_1 \geq \dots \geq \alpha_n)$ and $\gamma =
(\gamma_1 \geq \dots \geq \gamma_n)$ be weakly decreasing sequences of real
numbers.  There exist Hermitian $n\times n$ matrices $A$ and $C$ with
these eigenvalues such that $C-A$ is positive semidefinite of rank at
most one, if and only if $\gamma_1 \geq \alpha_1 \geq \gamma_2 \geq
\alpha_2 \geq \dots \geq \gamma_n \geq \alpha_n$.
\end{lemma}
\begin{proof}
  Set $\beta = (\beta_1,0,\dots,0)$ where $\beta_1 = \sum \gamma_i -
  \sum \alpha_i$, and assume that $\beta_1 \geq 0$.  We must show that
  there are Hermitian matrices $A$, $B$, and $C$ with eigenvalues
  $\alpha$, $\beta$, and $\gamma$ such that $A+B=C$ if and only if
  $\gamma_1 \geq \alpha_1 \geq \dots \geq \gamma_n \geq \alpha_n$.
  
  By approximating the eigenvalues with rational numbers and clearing
  denominators, we may assume that $\alpha$, $\beta$, and $\gamma$ are
  partitions.  In this case it follows from the work of Klyachko
  \cite{klyachko:stable} and Knutson and Tao
  \cite{knutson.tao:honeycomb} that the matrices $A,B,C$ exist
  precisely when the Littlewood-Richardson coefficient
  $c^\gamma_{\alpha \beta}$ is non-zero (see
  \cite[Thm.~11]{fulton:eigenvalues*1}).  This is equivalent to the
  specified inequalities by the Pieri rule.
\end{proof}

We remark that Lemma~\ref{L:rankone} can also be deduced directly from
Klyachko's result \cite[Thm.~1.2]{klyachko:stable} in a few more
lines.  Finally, we need the following statement, which is equivalent
to the Claim proved in \cite[p.~30]{fulton:eigenvalues*2}.

\begin{lemma}[Fulton] \label{L:fulton}
  Let $\alpha(1),\dots,\alpha(m)$ be weakly decreasing $n$-tuples of
  real numbers which satisfy \major{n}.  Suppose that for some
  sequence $(I(1),\dots,I(m)) \in S^n_t(m)$ we have $\sum_{s=1}^m
  \sum_{i \in I(s)} \alpha_i(s) = 0$.  For $1 \leq s \leq m$ we let
  $\alpha'(s)$ be the sequence of $\alpha_i(s)$ for $i \in I(s)$ and
  let $\alpha''(s)$ be the sequence of $\alpha_i(s)$ for $i \not \in
  I(s)$, both in weakly decreasing order.  Then $\{\alpha'(s)\}$
  satisfy \major{t} and $\{\alpha''(s)\}$ satisfy \major{n-t}.
\end{lemma}

We prove that the inequalities \major{n} and \rkcond{n}{r} are
sufficient by a `lexicographic' induction on $(n,r)$.  As the starting
point we take the cases where $r=0$, which are already known
\cite{klyachko:stable, belkale:local, knutson.tao.ea:honeycomb},
\cite[Thm.~17]{fulton:eigenvalues*1}.  For the induction step we let
$1 \leq r \leq n$ be given and assume that the inequalities are
sufficient in all cases where $n$ is smaller, as well as the cases
with the same $n$ and smaller $r$.  Using this hypothesis, we start by
proving the following fact.  Given two decreasing $n$-tuples $\alpha$
and $\beta$, we write $\alpha \geq \beta$ if $\alpha_i \geq \beta_i$
for all $i$.

\begin{lemma} \label{L:between}
  Let $\beta$, $\gamma$, and $\alpha(2),\dots,\alpha(m)$ be weakly
  decreasing $n$-tuples with $\beta \geq \gamma$, such that
  $\beta,\alpha(2),\dots,\alpha(m)$ satisfy \major{n} and
  $\gamma,\alpha(2),\dots,\alpha(m)$ satisfy \rkcond{n}{r}.  There
  exists a decreasing $n$-tuple $\alpha(1)$ such that $\beta \geq
  \alpha(1) \geq \gamma$ and $\alpha(1),\dots,\alpha(m)$ satisfy both
  \major{n} and \rkcond{n}{r}.
\end{lemma}
\begin{proof}
  We start by decreasing some entries of $\beta$ in the following way.
  First decrease $\beta_n$ until an inequality \major{n} becomes an
  equality, or until $\beta_n = \gamma_n$.  If the latter happens,
  then we continue by decreasing $\beta_{n-1}$ until an inequality
  \major{n} becomes an equality, or until $\beta_{n-1} =
  \gamma_{n-1}$.  If the latter happens we continue by decreasing
  $\beta_{n-2}$, etc.  If we are able to decrease all entries in
  $\beta$ so that $\beta = \gamma$, then we can use $\alpha(1) =
  \gamma$.
  
  Otherwise we may assume that for some sequence $(I(1),\dots,I(m))
  \in R^n_t(m)$ we have an equality $\sum_{i \in I(1)} \beta_i +
  \sum_{s=2}^m \sum_{i\in I(s)} \alpha_i(s) = 0$.  For each $s\geq 2$
  we let $\alpha'(s)$ be the decreasing $t$-tuple of numbers
  $\alpha_i(s)$ for $i \in I(s)$, and we let $\alpha''(s)$ be the
  decreasing $(n-t)$-tuple of numbers $\alpha_i(s)$ for $i \not \in
  I(s)$.  Similarly we define decreasing tuples
  $\beta'=(\beta_i)_{i\in I(1)}$, $\beta''=(\beta_i)_{i \not \in
    I(1)}$, and $\gamma''=(\gamma_i)_{i \not \in I(1)}$.  By
  Lemma~\ref{L:fulton} we know that $\beta',
  \alpha'(2),\dots,\alpha'(m)$ satisfy \major{t} and that $\beta'',
  \alpha''(2), \dots, \alpha''(m)$ satisfy \major{n-t}.  In
  particular, since the entries of the $t$-tuples add up to zero, we
  can find Hermitian $t \times t$ matrices $A'(1),\dots,A'(m)$ with
  eigenvalues $\gamma',\alpha'(2),\dots,\alpha'(m)$ such that $\sum
  A'(s) = 0$.
  
  We claim that the $(n-t)$-tuples $\gamma'', \alpha''(2), \dots,
  \alpha''(m)$ satisfy \rkcond{n-t}{r}.  This is clear if $n-t \leq
  r$.  Otherwise set $J(s) = \{ n+1-i \mid i \not \in I(s) \}$.  Since
  $\lambda(J(s))$ is the conjugate partition of $\lambda(I(s))$, it
  follows that $(J(1),\dots,J(m)) \in R^n_{n-t}(m)$.  For any sequence
  $(P(1),\dots,P(m)) \in R^{n-t-r}_x(m)$, we obtain from
  Lemma~\ref{L:index} that the sequence
  $(J(1)_{P(1)},\dots,J(m)_{P(m)})$ belongs to $S^{n-r}_x(m)$.  Notice
  that if $J(s) = \{ j_1 < j_2 < \dots < j_{n-t} \}$, then
  $\alpha''_{n-t+1-p}(s) = \alpha_{n+1-j_p}(s)$.  The claim therefore
  follows because
\begin{multline*} 
  \sum_{p\in P(1)} \gamma''_{n-t+1-p} + \sum_{s=2}^m \sum_{p\in P(s)}
  \alpha''_{n-t+1-p}(s) \ \ = \\ \sum_{j \in J(1)_{P(1)}} \gamma_{n+1-j} +
  \sum_{s=2}^m \sum_{j\in J(s)_{P(s)}} \alpha_{n+1-j}(s) \ \ \leq \ \ 0 \,.
\end{multline*}

By induction on $n$ there exists a decreasing $(n-t)$-tuple
$\alpha''(1)$ such that $\beta'' \geq \alpha''(1) \geq \gamma''$ and
$\alpha''(1),\dots,\alpha''(m)$ satisfy both of \major{n-t} and
\rkcond{n-t}{r}.  By the cases of Theorem~\ref{T:majorrank} that we
assume are true by induction, we can find Hermitian $(n-t)\times
(n-t)$ matrices $A''(1),\dots,A''(m)$ with eigenvalues
$\alpha''(1),\dots,\alpha''(m)$ and with positive semidefinite sum of
rank at most $r$.  We can finally take $\alpha(1)$ to be the
eigenvalues of $A'(1) \oplus A''(1)$.
\end{proof}

We can now finish the proof that the inequalities of
Theorem~\ref{T:majorrank} are sufficient.  Let $\gamma =
(\alpha_2(1),\alpha_3(1),\dots,\alpha_n(1),M)$ for some large
negative number $M \ll 0$.  We claim that when $M$ is sufficiently
small, the $n$-tuples $\gamma,\alpha(2),\dots,\alpha(m)$ satisfy
\rkcond{n}{r-1}.  In fact, let $(P(1),\dots,P(m)) \in R^{n-r+1}_t(m)$.
If $1 \in P(1)$ then the inequality for this sequence holds by choise
of $M$.  Otherwise we have that $(Q,P(2),\dots,P(m)) \in R^{n-r}_t(m)$
where $Q = \{p-1 \mid p \in P(1)\}$, and the required inequality
follows because
\[ \sum_{q \in Q} \alpha_{n+1-q}(1) + \sum_{s=2}^m \sum_{p \in P(s)}
   \alpha_{n+1-p}(s) \ \leq \ 0 \,.
\]

By Lemma~\ref{L:between} we may now find a decreasing $n$-tuple
$\tilde \alpha(1)$ with $\alpha(1) \geq \tilde \alpha(1) \geq \gamma$,
such that $\tilde \alpha(1), \alpha(2), \dots, \alpha(m)$ satisfy
\major{n} and \rkcond{n}{r-1}.  By induction on $r$ there exist
Hermitian $n \times n$ matrices $\tilde A(1), A(2), \dots, A(m)$ with
eigenvalues $\tilde \alpha(1),\linebreak \alpha(2), \dots, \alpha(m)$,
such that $\tilde A(1)+A(2)+\dots+A(m)$ is positive semidefinite of
rank at most $r-1$.  Finally, using Lemma~\ref{L:rankone} and the
choise of $\gamma$ we may find a Hermitian matrix $A(1)$ with
eigenvalues $\alpha(1)$ such that $A(1)-\tilde A(1)$ is positive
semidefinite of rank at most 1.  The matrices $A(1),A(2),\dots,A(m)$
now satisfy the requirements.

\section{Minimality of the inequalities}

In this section we prove that when $r \geq 1$ and $m \geq 3$, the
inequalities \decr{}, \major{n}, and \rkcond{n}{r} are independent,
thereby proving the last statement of Theorem~\ref{T:majorrank}.  It
is enough to show that for each inequality among \major{n} or
\rkcond{n}{r}, there exist strictly decreasing $n$-tuples
$\alpha(1),\dots,\alpha(m)$ such that the given inequality is an
equality and all other inequalities \major{n} and \rkcond{n}{r} are
strict.  In addition we must show that for each $1 \leq i \leq n-1$
there exist $\alpha(1) = (\alpha_1(1) > \dots > \alpha_i(1) =
\alpha_{i+1}(1) > \dots > \alpha_n(1))$ and strictly decreasing
$n$-tuples $\alpha(2),\dots,\alpha(m)$, such that all inequalities
\major{n} and \rkcond{n}{r} are strict.

We start with the latter case.  If $n=2$ we can take $\alpha(1) =
(0,0)$ and $\alpha(s) = (2,-1)$ for $2 \leq s \leq m$.  For $n \geq
3$, it was shown in \cite[Lemma 1]{fulton:eigenvalues} that the
$n$-tuples $\beta(1)=\beta(2)=\dots=\beta(m) =
(n-1,n-3,\dots,3-n,1-n)$ satisfy that $\sum_{s=1}^m \sum_{i\in I(s)}
\beta_i(s) \geq 2$ for all sequences $(I(1),\dots,I(m)) \in R^n_t(m)$
of subsets of cardinality $t < n$.  In fact, this follows because
$\sum_{s=1}^m \sum_{i\in I(s)} i = t(n-t) + m\binom{t+1}{2}$.  Using
this fact, one easily checks that both \major{n} and \rkcond{n}{r} are
strict for $\alpha(1) = (n-1,n-3,\dots,n-2i,n-2i,\dots,3-n,1-n)$, with
$n-2i$ as the $i$th and $i+1$st entries, and $\alpha(2) = \dots =
\alpha(m) = (n,n-3,n-5,\dots,3-n,1-n)$.

Now consider an inequality from \major{n}, given by a sequence
$(I(1),\dots,I(m)) \in R^n_t(m)$.  By
\cite[Thm.~9]{knutson.tao.ea:honeycomb} we can choose strictly
decreasing $n$-tuples $\alpha(1),\dots,\alpha(m)$ such that
$\sum_{s=1}^m \sum_{i=1}^n \alpha_i(s) = \sum_{s=1}^m \sum_{i \in
  I(s)} \alpha_i(s) = 0$ and all other inequalities \major{n} are
strict.  If $(P(1),\dots,P(m)) \in R^{n-r}_x(m)$ then we have
$(Q,P(2),\dots,P(m)) \in R^n_x(m)$ where $Q = \{p+r \mid p \in
P(1)\}$.  Since the negated $n$-tuples $\tilde \alpha(1),\dots, \tilde
\alpha(m)$ given by $\tilde \alpha(s) = (-\alpha_n(s) > \dots >
-\alpha_1(s))$ must satisfy \major{n}, we obtain that $\sum_{s=1}^m
\sum_{p \in P(s)} \alpha_{n+1-p}(s) < \sum_{q\in Q} \alpha_{n+1-q}(1)
+ \sum_{s=2}^m \sum_{p\in P(s)} \alpha_{n+1-p}(s) \leq 0$.  This shows
that the inequalities \rkcond{n}{r} are strict.  If $t < n$ we may
finally replace $\alpha_{i_0}(1)$ with $\alpha_{i_0}(1) + \epsilon$,
where $i_0 \not \in I(1)$, to obtain that $\sum_{s=1}^m \sum_{i=1}^n
\alpha_i(s) > 0$.

At last we consider an inequality of \rkcond{n}{r} given by a sequence
$(P(1),\dots,P(m)) \in R^{n-r}_x(m)$.  We once more apply
\cite[Thm.~9]{knutson.tao.ea:honeycomb} to obtain strictly decreasing
$(n-r)$-tuples $\beta(1),\dots,\beta(m)$ such that $\sum_{s=1}^m
\sum_{p=1}^{n-r} \beta_p(s) = \sum_{s=1}^m \sum_{p\in P(s)}
\beta_{p}(s) = 0$, and all other inequalities of \major{n-r} are
strict.  Set $\alpha(s) = (N+r, N+r-1,\linebreak \dots, N+1,
-\beta_{n-r}(s), \dots, -\beta_1(s))$ for $1 \leq s \leq m$, where $N
\gg 0$ is a large number.  Then the $n$-tuples
$\alpha(1),\dots,\alpha(m)$ strictly satisfy all inequalities from
\rkcond{n}{r}, except for the equalities $\sum_{s=1}^m
\sum_{p=1}^{n-r} \alpha_{n+1-p}(s) = \sum_{s=1}^m \sum_{p \in P(s)}
\alpha_{n+1-p}(s) = 0$.  We must show that $\sum_{s=1}^m \sum_{i \in
  I(s)} \alpha_i(s) > 0$ for every sequence $(I(1),\dots,I(m)) \in
R^n_t(m)$.  If $I(1) \cap [r] \neq \emptyset$ then this follows from
our choise of $N$.  Otherwise we have $(J,I(2),\dots,I(m)) \in
R^{n-r}_t(m)$ where $J = \{ i-r \mid i \in I(1) \}$.  Since
$\alpha_i(s) > -\beta_{n-r+1-i}(s)$ for $i \in [n-r]$, we obtain that
$\sum_{s=1}^m \sum_{i\in I(s)} \alpha_i(s) > \sum_{i\in J}
(-\beta_{n-r+1-i}(1)) + \sum_{s=2}^m \sum_{i\in I(s)}
(-\beta_{n-r+1-i}(s)) \geq 0$.  Finally, if $x \neq n-r$ we replace
$\alpha_{n+1-p_0}(1)$ with $\alpha_{n+1-p_0}(1) - \epsilon$, $p_0 \not
\in P(1)$, to obtain a strict inequality $\sum_{s=1}^m
\sum_{p=1}^{n-r} \alpha_{n+1-p}(s) < 0$.  This completes the proof
that the inequalities are independent.


\begin{thebibliography}{1}

\bibitem{belkale:local}
P.~Belkale, \emph{Local systems on ${\mathbb P}^1 \smallsetminus s$ for $s$ a
  finite set}, Ph.D. thesis, University of Chicago, 1999.

\bibitem{friedland:finite}
S.~Friedland, \emph{Finite and infinite dimensional generalizations of
  {K}lyachko's theorem}, Linear Algebra Appl. \textbf{319} (2000), no.~1-3,
  3--22. \MR{2002a:15023}

\bibitem{fulton:eigenvalues}
W.~Fulton, \emph{Eigenvalues of sums of {H}ermitian matrices (after {A}.
  {K}lyachko)}, Ast\'erisque (1998), no.~252, Exp.\ No.\ 845, 5, 255--269,
  S\'eminaire Bourbaki. Vol.\ 1997/98. \MR{2000e:14092}

\bibitem{fulton:eigenvalues*1}
\bysame, \emph{Eigenvalues, invariant factors, highest weights, and {S}chubert
  calculus}, Bull. Amer. Math. Soc. (N.S.) \textbf{37} (2000), no.~3, 209--249
  (electronic). \MR{2001g:15023}

\bibitem{fulton:eigenvalues*2}
\bysame, \emph{Eigenvalues of majorized {H}ermitian matrices and
  {L}ittlewood-{R}ichardson coefficients}, Linear Algebra Appl. \textbf{319}
  (2000), no.~1-3, 23--36.  \MR{2002a:15024}

\bibitem{klyachko:stable}
A.~A. Klyachko, \emph{Stable bundles, representation theory and {H}ermitian
  operators}, Selecta Math. (N.S.) \textbf{4} (1998), no.~3, 419--445.
  \MR{2000b:14054}

\bibitem{knutson.tao:honeycomb}
A.~Knutson and T.~Tao, \emph{The honeycomb model of 
  {${\rm GL}\sb n({\mathbb C})$} tensor products. 
  {I}. {P}roof of the saturation conjecture}, J. Amer. Math.
  Soc. \textbf{12} (1999), no.~4, 1055--1090. \MR{2000c:20066}

\bibitem{knutson.tao.ea:honeycomb}
A.~Knutson, T.~Tao, and C.~Woodward, \emph{The honeycomb model of {${\rm GL}\sb
  n(\mathbb C)$} tensor products. {II}. {P}uzzles determine facets of the
  {L}ittlewood-{R}ichardson cone}, J. Amer. Math. Soc. \textbf{17} (2004),
  no.~1, 19--48 (electronic). \MR{2 015 329}

\bibitem{weyl:das}
H.~Weyl, \emph{Das asymptotische {V}erteilungsgesetz der {E}igenwerte lineare
  partieller {D}ifferentialgleichungen}, Math. Ann. \textbf{71} (1912),
  441--479.

\end{thebibliography}

\providecommand{\bysame}{\leavevmode\hbox to3em{\hrulefill}\thinspace}
\providecommand{\MR}{\relax\ifhmode\unskip\space\fi MR }
\providecommand{\MRhref}[2]{%
  \href{http://www.ams.org/mathscinet-getitem?mr=#1}{#2}
}
\providecommand{\href}[2]{#2}


\end{document}